\documentclass{article}
\usepackage{graphicx} 
\usepackage{amsmath}
\usepackage[portrait, margin=3cm]{geometry}

\usepackage[square,numbers]{natbib}
\newcommand{\be}{\ensuremath{\mathbf{e}}}

\newcommand{\beef}{\ensuremath{\mathbf{f}}}

\newcommand{\beem}{\ensuremath{\mathbf{m}}}

\newcommand{\bu}{\ensuremath{\mathbf{u}}}
\newcommand{\bv}{\ensuremath{\mathbf{v}}}

\newcommand{\bF}{\ensuremath{\mathbf{F}}}

\newcommand{\bM}{\ensuremath{\mathbf{M}}}

\newcommand{\bW}{\ensuremath{\mathbf{W}}}
\newcommand{\bX}{\ensuremath{\mathbf{X}}}

\newcommand{\balpha}{\ensuremath{\boldsymbol{\alpha}}}
\newcommand{\bbeta}{\ensuremath{\boldsymbol{\beta}}}

\newcommand{\btheta}{\ensuremath{\boldsymbol{\theta}}}
\newcommand{\bTheta}{\ensuremath{\boldsymbol{\Theta}}}

\newcommand{\etal}{{\it{et al}}.}

\newcommand{\bigR}{\ensuremath{\mathbf{R}}}
\newcommand{\bigF}{\ensuremath{\mathbf{F}}}
\newcommand{\bigM}{\ensuremath{\mathbf{M}}}

\newcommand{\bigY}{\ensuremath{\mathbf{Y}}}
\newcommand{\bigU}{\ensuremath{\mathbf{U}}}
\newcommand{\bigW}{\ensuremath{\mathbf{W}}}

\newcommand{\bigTheta}{\Theta}

\usepackage{bm}

\title{A Complete Graphic Statics for Rigid-Jointed 3D Frames. \\Part 3: Loops for Kinematics}
\author{Allan McRobie, fam20@cam.ac.uk \\
Dept. of Engineering, Cambridge University, CB2 1PZ\\
https://orcid.org/0000-0002-6610-5927}

\begin{document}

\maketitle

\section*{\begin{center} Abstract \end{center}}
In Part 3 of this sequence of papers, the kinematic behaviour of 3D frame structures is described using the loop formalism that was developed in Part 2 to describe equilibrium. There, the notions of polygons, polyhedra and polytopes that form the geometric toolbox underlying graphic statics were replaced by the more general concept of CW-complexes from algebraic homology. The six components of the stress resultant acting on any cut face of a bar in a rigidly-jointed framework were represented by the oriented bivector areas of the six projections of a loop in a 4D-space, with three components representing the force and three components representing the moment.  In this paper, projected areas of loops in 4D will represent kinematic variables, with three projected areas representing the displacement of a point on the frame, and three other projected areas representing the rotation of the structure at that point. 

The 4D setting for the theory consists of the usual three dimensions of physical space together with a fourth dimension for the stress function.
Virtual Work then manifests as a top form (an oriented 4-volume) in this 4D setting, being  the integral over the structure of the wedge product of bivectors representing the local equilibrium and kinematic variables.

\section{Introduction}
Previous parts of this sequence of papers~\cite{McRobieArXiv1, McRobieArXiv2} have been concerned with the equilibrium of forces and moments  in structural frames. Here we consider the kinematics of such structures, extending the 
McRobie \etal~\cite{McRobieRSOS2} description of the kinematics of pin-jointed trusses to the case of moment-resisting frames. 

\begin{figure} [t!] \centering
\includegraphics[width = 0.88\textwidth]{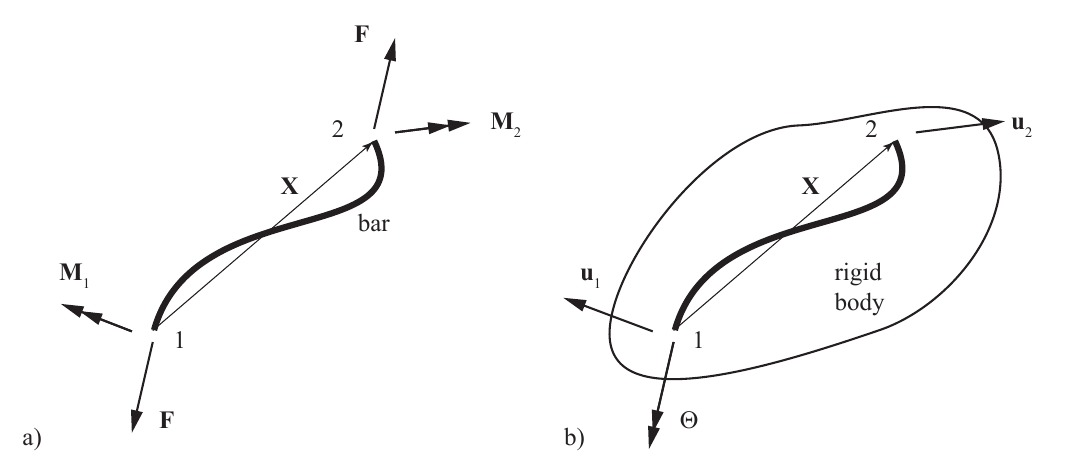}
\caption{a) A bar subject to end forces and moments.
b) A rigid body motion.}
\label{BasicBar}
\end{figure}

The starting point is the recognition of the obvious similarity between the equations describing the equilibrium of a rigid body and its infinitesimal motion.
For example,   equilibrium of a bar subject to end forces $\bF$ and moments $\bM_1$ and $\bM_2$ (see Fig.~\ref{BasicBar}), requires
\begin{equation}
\bM_2 = \bM_1 + \bX \times \bF
\end{equation}
Similarly, for a rigid body undergoing an infinitesimal displacement involving a rotation $\btheta$, the displacements $\bu_1, \bu_2$ at any two points may be written (again in traditional vector notation) as
\begin{equation}
\bu_2 = \bu_1 + \bX \times \btheta
\end{equation}

Even though the frames that we shall consider may be more generally flexible bodies, it is clear that there is some
similarity in the way that moments and displacements are treated, and likewise for forces and rotations.  This may be somewhat counterintuitive at first, when there may be a suspicion that perhaps rotations and moments are somehow similar, and that displacements and forces somehow resemble each other. We now take these ideas from the 3D description to the 4D Legendre transform description developed in the earlier papers here ~\cite{McRobieArXiv1,McRobieArXiv2}, where an additional dimension has been added to carry various forms of stress function.

\subsection{Notation}
\begin{figure} [b!] \centering
\includegraphics[width = 0.88\textwidth]{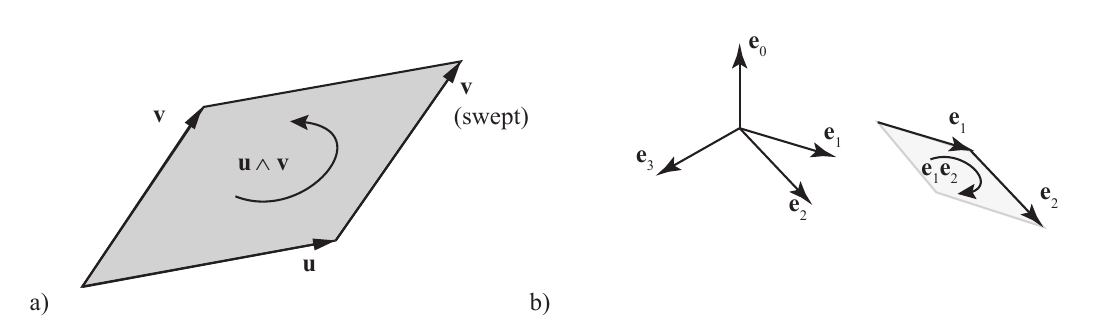}
\caption{a) The wedge product of general vectors $\bu$ and $\bv$. b) The unit bivector $\be_1\be_2$
}
\label{bivector}
\end{figure}

Since the vector cross product is not defined in 4D we use instead the wedge product of exterior algebra, which works in any dimension.  The wedge product of vectors $\bu$ and $\bv$ is the bivector  $\bu \wedge \bv$, which corresponds to the oriented area of the parallelogram swept out when the vector $\bv$ is swept along $\bu$ (see Fig.~\ref{bivector}a). The orientation is defined by traversing the boundary first along side $\bu$ and then along the side created by the swept $\bv$. In both the body space and the stress space, we  shall choose an orthonormal set of basis vectors 
$\be_0$, $\be_1$, $\be_2$, $\be_3$, with the $\be_0$ direction associated with the stress functions $F$ and $\Phi$ of the Legendre transformation between the two spaces.

For notational convenience we may write the wedge product of basis vectors $\be_i$ and $\be_j$ as the Clifford product $\be_i \be_j$ whenever $i$ and $j$ are different. This follows from the Clifford product definition $\be_i\be_j = \be_i . \be_j + \be_i \wedge \be_j$ with the dot product being zero by orthonormality.
We make no use of the Clifford product other than to save writing wedge symbols.

It is difficult to draw four dimensional objects. Nevertheless unit bivectors and their orientations can be readily represented diagrammatically, as shown in Fig.~\ref{bivector}b for $\be_1\be_2$.

\subsection{Extension of 4D loop formalism to kinematics}

\begin{figure} [t!] \centering
\includegraphics[width = \textwidth]{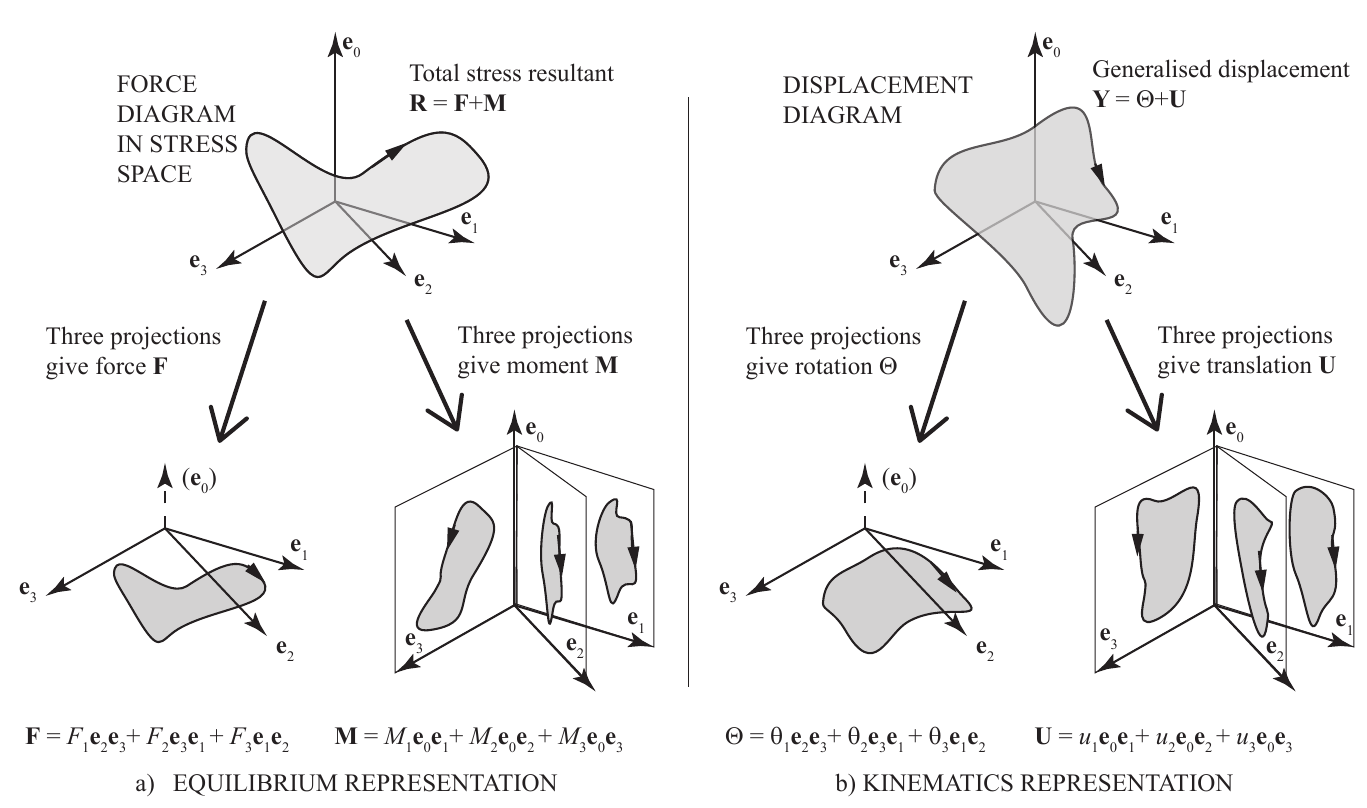}
\caption{The loop representation for: a) equilibrium; b) kinematics.}
\label{EquilibRep}
\end{figure}

In the 4D description of the equilibrium of frame structures, a general stress resultant $\bigR$ (i.e. the forces and the moments) at any cut face of a bar is represented by the six independent orthogonal components of the oriented bivector area of a general loop or set of loops in the 4D stress space (see Fig.~\ref{EquilibRep}a). That is, we write
\begin{eqnarray}
\bigR & = & \bigF + \bigM = \bigF + \be_0 \beem \\
\mathrm{where} \ \ \  \bigF & = & F_1 \be_2 \be_3 + F_2 \be_3 \be_1 + F_3 \be_1 \be_2 \\
\mathrm{and}   \ \ \  \bigM & = & M_1 \be_0 \be_1 + M_2 \be_0 \be_2 + M_3 \be_0 \be_3
\end{eqnarray}
Here, $\bigR$, $\bigF$ and $\bigM$ are bivectors. We may later wish to make use of the (possibly more familiar) 3-vectors $\beem = M_1\be_1 + M_2\be_2 + M_3\be_3$ and $\beef = F_1\be_1 + F_2\be_2 + F_3\be_3$.

An analogous kinematic description may be obtained by writing the generalised infinitesimal displacement $\bigY$ of any point on the frame as the six bivector components of a general loop or set of loops, using 
\begin{eqnarray}
\bigY & = & \bigTheta + \bigU = \bigTheta + \be_0 \bu \\
\mathrm{where} \ \ \  \bigTheta & = & \theta_1 \be_2 \be_3 + \theta_2 \be_3 \be_1 + \theta_3 \be_1 \be_2 \\
\mathrm{and}   \ \ \  \bigU & = & u_1 \be_0 \be_1 + u_2 \be_0 \be_2 + u_3 \be_0 \be_3
\end{eqnarray}

Similar to the equilibrium description, $\bigY$, $\bigTheta$ and $\bigU$ are bivectors. Again. the more familiar 3-vectors would be $\bu = u_1\be_1 + u_2\be_2 + u_3\be_3$ and $\btheta = \theta_1\be_1 + \theta_2\be_2 + \theta_3\be_3$.
These objects lie in the 4D displacement space, which may be identified with the body space.

This all accords with Felix Klein's observation that ``All considerations that apply to forces that act on a rigid body
can be applied in a completely analogous form to the infinitely small rotations that such a
rigid body performs, and conversely''~\cite{Klein2}. Here, though, we are extending the idea to non-rigid bodies, and will make a similar analogy between moments and translations. As already noted, this is the converse of any seemingly more natural temptation to pair moments with rotations, and forces with displacements.

\begin{figure} [hb!] \centering
\includegraphics[width = 0.88\textwidth]{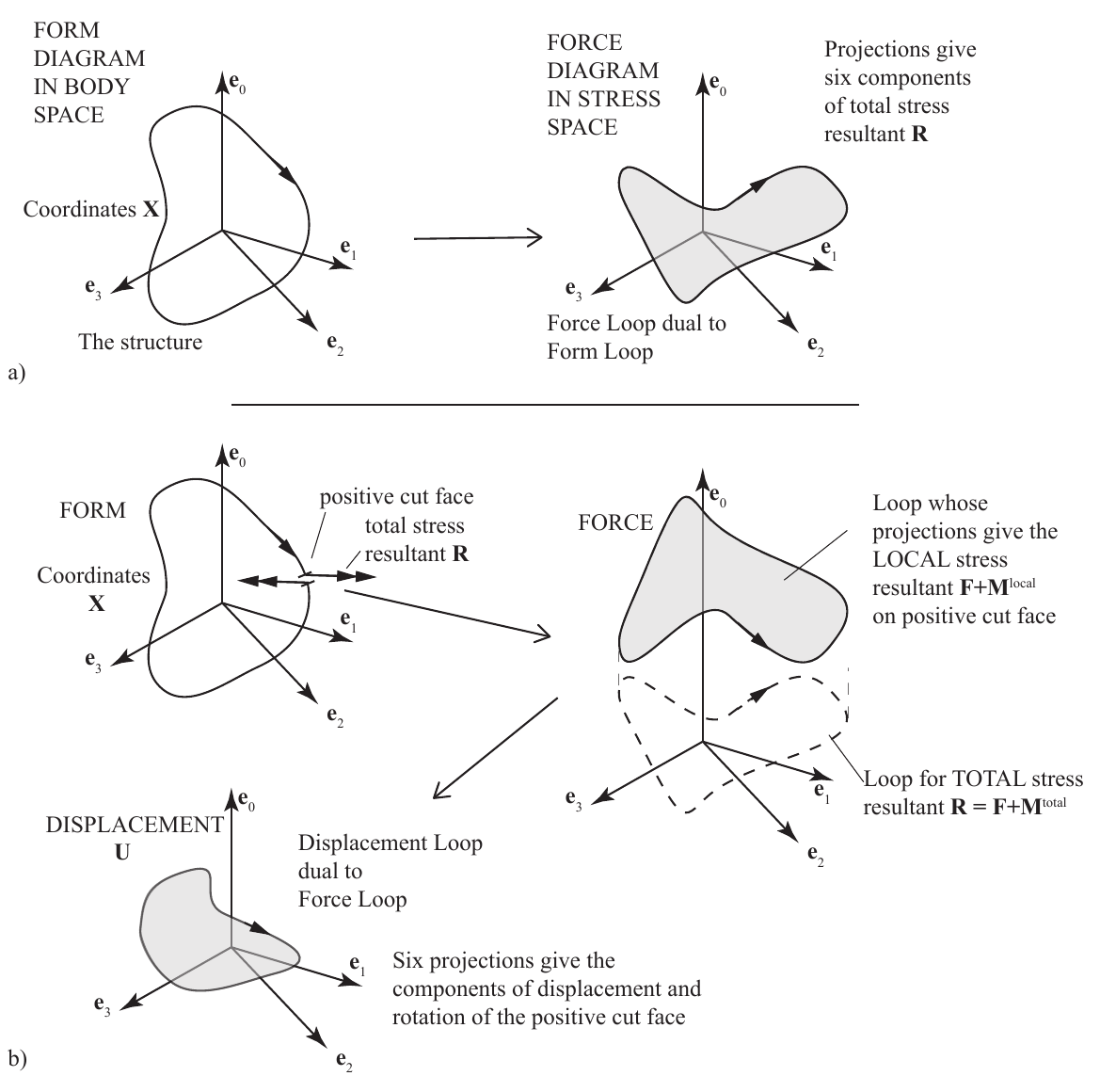}
\caption{a) A loop of a simple moment-resisting frame structure and its dual force loop whose six projections give the components of the total stress resultant $\bigR$ acting on any cut face of the structural loop. The total moment includes the contribution from the force acting at a lever arm about the origin. By elementary equilibrium, this is the same at all positive cut faces.  
b) Given a positive cut face on the structural loop, an adjusted force loop can be defined whose projections give the local stress resultant where the moments are the local bending and torsional moments about the cut face. These moments vary around the structural loop.  Dual to any such local force loop a generalised displacement loop $\bigU$ can be defined whose six projections give the rotation and displacement of the  structure at the cut.}
\label{FormForceDisp}
\end{figure}

The various geometric objects are illustrated in Fig.~\ref{FormForceDisp}. The frame structure is decomposed into a set of oriented loops, as described in McRobie~\cite{McRobieArXiv1}, and one such loop is shown in Fig.~\ref{FormForceDisp}a. Dual to this form loop there is a force loop in the stress space which defines the state of stress on any positive cut face of the form loop. The six independent projections of the force loop define the six components of the stress resultant at any positive cut face. The moments so represented are the total moments, which include the contribution of the forces and their lever arms about the origin. Since there are no forces applied externally to the bar, it follows from elementary equilibrium that the total moment is the same at any positive cut face on the form loop, such that there is just a single force loop dual to each form loop. It is explained in McRobie~\cite{McRobieArXiv1,McRobieArXiv2} how the local bending moments and torsions may be represented by the projections of a continuum of related loops in the stress space. These are obtained by plotting the original Maxwell-Rankine stress function $F$ in the $\be_0$ direction above the stress space $\xi, \eta, \zeta$, rather than the Legendre dual stress function $\Phi$. There is such an adjusted force loop, as illustrated in Fig.~\ref{FormForceDisp}b, corresponding to any cut face of the form loop. Unlike the unique loop which represents force and total moment at any cut, there is thus a continuum of adjusted force loops which capture how the bending and torsional moments vary around the loop.

Dual to each adjusted force loop there is a displacement loop whose projections give the translation and rotation of the cut face on the structural loop. 
One such displacement loop is shown in Fig.~\ref{FormForceDisp}b.
The structure is not a rigid body and thus, similar to the adjusted force loops, there is an infinite set of displacement loops, one for each possible cut through the form loop. These loops do not necessarily form a continuum: it is possible to have discontinuities in displacement and rotation. Nevertheless, there is such a loop for each cut.

We now have all the elements required for a geometric treatment of Virtual Work. The Principle of Virtual Work is arguably the fundamental object in structural analysis, and we now illustrate how it is manifested in the loop formalism.

\section{Virtual Work}
An oriented loop $\bigR$ in the 4D stress space represents the six components of the general stress resultant of force and moment at any cut on the structural loop. We may also associate an oriented loop $\bigY$ in the body space to represent the six independent components of the generalised infinitesimal translation and rotation. The Virtual Work $\bigW$ associated with the stress resultants undergoing that infinitesimal motion is then
\begin{eqnarray}
\bigW & = & \bigR \wedge \bigY  = (\bigF + \bigM) \wedge (\bigTheta + \bigU)  \nonumber \\
      & = &  W \be_0\be_1\be_2\be_3 \nonumber \\
 \mathrm{where} \quad W     & = & \beef.\bu + \beem.\btheta \label{eqn:W}
\end{eqnarray}
That is, the Virtual Work is a top form, the oriented volume of a four dimensional region. The magnitude of this volume is given by the
sum of the 3-vector dot products $\beef.\bu + \beem.\btheta$. This is a familiar quantity in traditional structural mechanics, but here
it has been given a new realisation as an oriented 4-volume. There is a rather pleasing neatness to the way that the cross terms $\bigF \wedge \bigTheta$  and $\bigM \wedge \bigU$ disappear, and the way that the remaining terms each deliver a 4-volume $\be_0\be_1\be_2\be_3$ whose magnitude is given by the familiar vector dot products. The explicit working is as follows: the wedge product $\balpha \wedge \bbeta$ of multivectors $\balpha = \be_{i_1} \wedge \ldots  \wedge \be_{i_p}$ and $\bbeta = \be_{j_1} \wedge \ldots  \wedge \be_{j_q}$
is given by
\begin{equation*}
    \balpha \wedge \bbeta =  \be_{i_1} \wedge \ldots  \wedge \be_{i_p} \wedge  \be_{j_1} \wedge \ldots  \wedge \be_{j_q}
\end{equation*}
if the indices $i_1, \ldots i_{p}, j_1, \ldots j_{q}$ are distinct, and equals zero otherwise (due to the presence of terms, after re-ordering,  of the form $\be_k \wedge \be_k = 0 $). It follows that the wedge product of the generalised force component $F_1 \be_2 \be_3$ is zero with all components of the generalised  displacement $\bigY$ except for the component $u_1 \be_0 \be_1$.  This gives the contribution
$F_1 u_1 \be_0\be_1\be_2\be_3$ to the Virtual Work. The other generalised force components similarly pick out just their conjugate components of the generalised displacement, leading to equation~(\ref{eqn:W}).  Since the four basis vectors are distinct and orthogonal, equation~(\ref{eqn:W}) uses the Clifford product notation $\be_i\be_j$ as shorthand for the wedge product $\be_i \wedge \be_j$.

\subsection{Relation to the Minkowski sum of truss diagrams}
The Minkowski sum~\cite{McRobieMink} has proven to be a useful visualisation tool for illustrating a range of important topics in graphic statics. In its simplest form, it allows the form and force diagrams of 2D trusses to be combined in a way that allows the load path to be visualised. The polygonal elements of each diagram are conjoined by a set of intervening rectangles. Since the rectangle side lengths are tension $T$ and length $L$, the area represents the contribution $T.L$ of that bar to the structure's load path. It follows trivially from the geometry of the construction that the total oriented area of all rectangles must be zero, and this corresponds to Maxwell's Load Path Theorem~\cite{maxwell1870a}. The concept generalises readily to 3D trusses and their Rankine reciprocals~\cite{McRobieMink}, with the Load Path now represented by right prisms of length $L$ and cross-sectional area $T$ separating the polyhedral cells of the form diagram from those of the Rankine reciprocal force diagram. 

In McRobie \etal~\cite{McRobieRSOS2} the idea was extended to allow Virtual Work to be visualised for 2D and 3D trusses. A new object was defined, being the assemblage of the Minkowksi sum of corresponding components of the force diagram and the {\it vector displacement diagram} (rather than the form diagram). The graph of the displacement diagram has the same topology as the structural form diagram, but the nodal locations are given by the nodal displacement vectors $\bu$ rather than the structural coordinates $\bX$. Since nodal displacements are not orthogonal to bar tensions, the components of the displacement and force diagrams are separated, for the 2D case, by parallelograms rather than rectangles, and for the 3D case, by skew prisms. The area (in 2D) or volume (in 3D) of these skewed objects give that bar's contribution to the internal Virtual Work of the system. In this ansatz, there is no external virtual work as there are no external forces, any such having been replaced by a set  of additional members aligned along their lines of action, such that the external forces are now subsumed within a state of self-stress of a larger system. It follows trivially from the geometry of the construction that the sum of the oriented areas (in 2D) or oriented volumes (in 3D) that separate the two diagrams must be zero and this is a geometric representation of the Principle of Virtual Work for trusses.

Here we consider how this visual approach to Virtual Work may be manifested in the loop formalism. At first, the representation of rotations by the oriented areas of loops may seem an unusual starting point. However, as we shall describe, one possible approach is to consider the case of linear elastic material behaviour and to interpret the (scaled) bending moment diagram (BMD) as a loop. Scaled by the flexural rigidity, the BMD gives a graph of the curvature due to flexure. The area under this graph corresponds to the change of slope along the beam caused by the bending. That is, the beam end rotations are naturally related to the oriented area of a loop (the scaled BMD). This is the link between rotations and loop areas.

What follows is simply an alternative geometric description, using the loop formalism, of the Force Method which is a familiar mainstay of structural analysis, wherein the Principle of Virtual Work is used to dot a real displacement system with a virtual equilibrium system.

\subsection{ Virtual Work for beam flexure: A geometric statement of integration by parts (Dec 2018)}
The expression of Virtual Work that we shall develop will contain terms of the form
\begin{equation}
\int M \kappa \ ds
\end{equation}
where $M$ is some bending moment from an equilibrium system and $\kappa$ is some conjugate curvature from a compatibility system. We illustrate how this is replicated in the new paradigm with some simple examples, starting with the 2D case of an Euler-Bernoulli beams undergoing small lateral deformations.

\subsubsection{The structure and the equilibrium system}
\begin{figure} [ht!] \centering
\includegraphics[width = \textwidth]{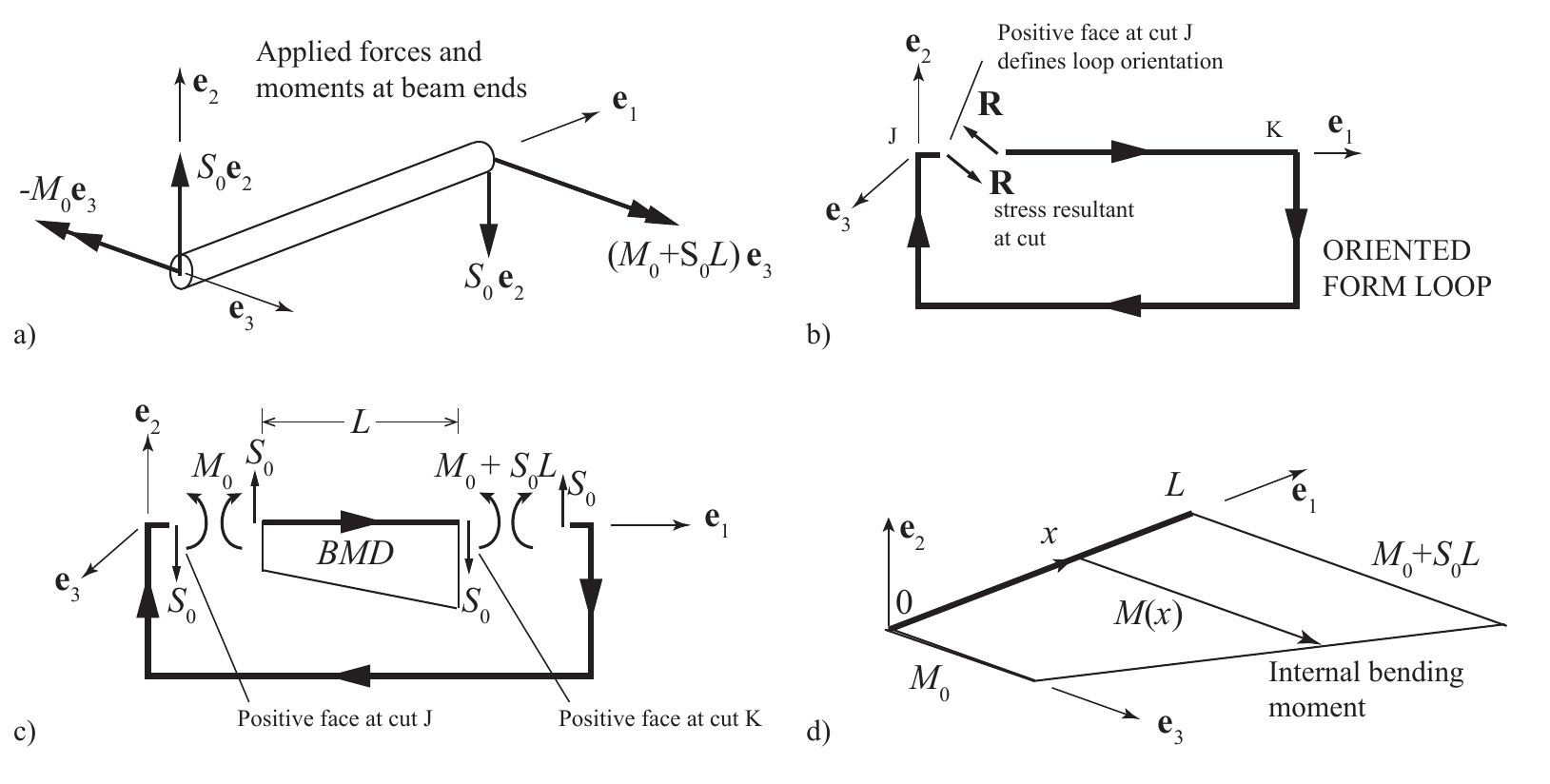}
\caption{a) A beam subject to applied end shears and end moments. b) The form diagram is an oriented loop created by connecting the beam ends with an arbitrary return path. The choice of positive face at cut J defines the loop orientation. c) Traditionally the bending moment diagram (BMD) is drawn in the plane of bending, $\be_1\be_2$ here. d) The internal bending moment $M(x)$ varies along the beam, and is plotted as a vector in the direction $\be_3$ normal to the plane of the frame. This is perpendicular to the traditional BMD representation, but accords with the usual vector notation for a moment.}
\label{basicbeamloop}
\end{figure}

Consider a straight beam oriented along the $\be_1$ axis, parameterised by the coordinate $x$, from $x = 0$ at the start node $J$ to $x= L$ at the end node $K$.
This beam is subject to end moments and end shears, as shown in  Fig.~\ref{basicbeamloop}a. The loop representation of the form diagram is shown in Fig.~\ref{basicbeamloop}b, where an arbitrary return path has been added between the beam ends to create the closed loop. The location of the cut on this loop is chosen to be just along the beam at node $J$, which is also the origin of the coordinate system. The positive face is selected as that facing along $\be_1$, towards node $K$. This defines the orientation of the form loop. 

\begin{figure} [ht!] \centering
\includegraphics[width = 0.83\textwidth]{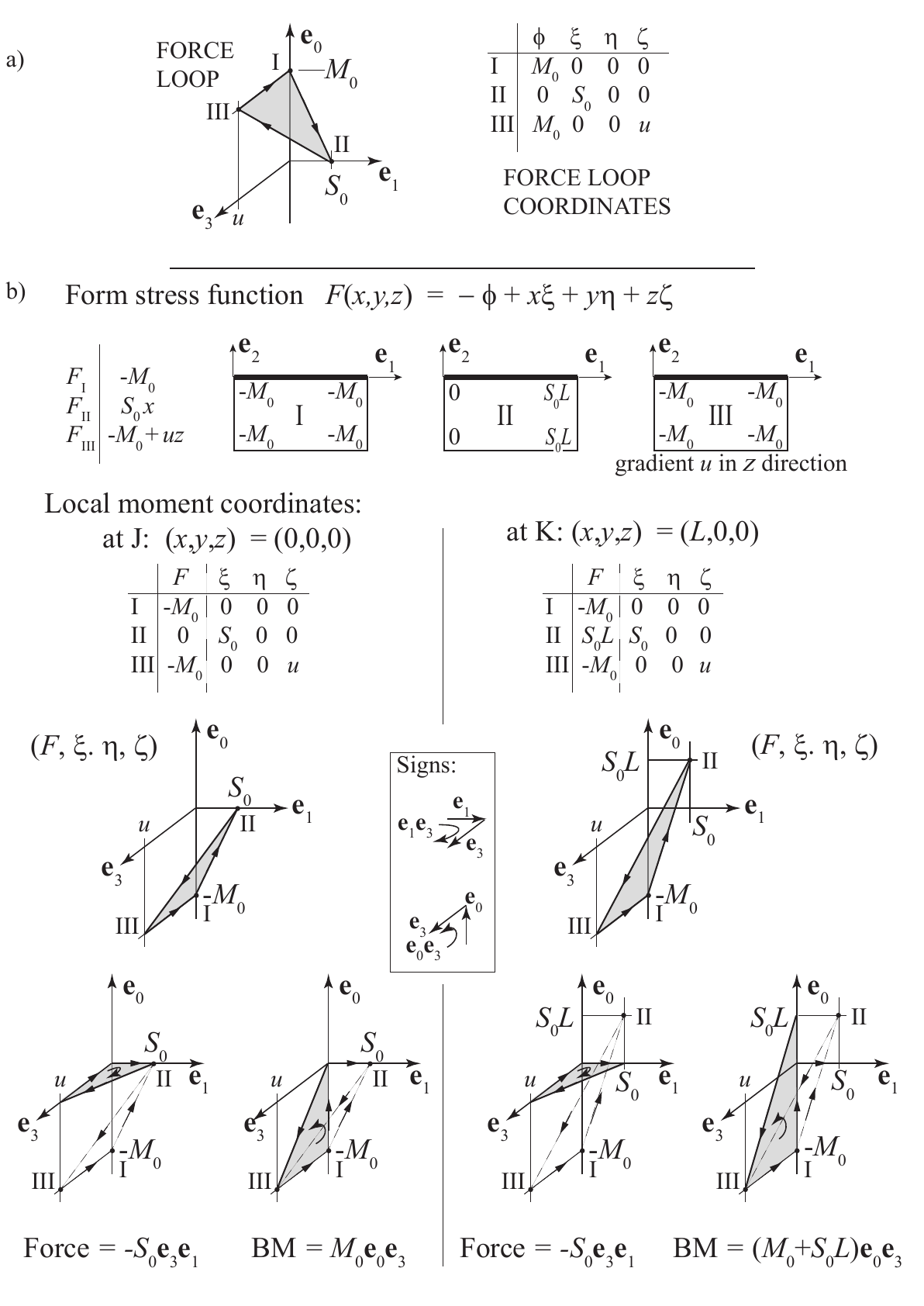}
\caption{a) The force diagram is any oriented loop in stress space $(\phi, \xi, \eta, \zeta)$ whose six projections give the total stress resultant at J.  Here we choose a simple triangular loop. The coordinate $\eta$ is zero for all points on the loop, thus the $\be_2$ dimension is suppressed, simplifying the diagram. b)  The local bending moments are given by the projected areas on the $\be_0\be_3$ plane of the loops in the hybrid $(F, \xi, \eta, \zeta)$ space. The appropriate projected area varies from $M_0$ at $J$ to $M_0+S_0L$ at $K$, as desired.}
\label{ForceDiagramforBeam}
\end{figure}

For the force diagram, we need to select a force loop which represents the total force and total moment about the origin of the stresses acting at the positive cut face at node $J$. Among the infinity of possibilities we select the simple triangular loop shown in Fig.~\ref{ForceDiagramforBeam}a. Its coordinates are defined by the values of $S_0, M_0$ and $u$, where $u$ is a quantity with units of force, and here we have chosen $u = 2$ units in order to simplify calculations by cancelling with the leading half in the formula for triangular areas.

As per McRobie~\cite{McRobieArXiv1}, local bending moments are given by projections of the hybrid loop $(F, \xi, \eta, \zeta)$ which uses the original stress function $F$ rather than the dual $\phi$. This is shown at length in Fig.~\ref{ForceDiagramforBeam}b, with areas projected on the $\be_0\be_3$ plane rising from $M_0$ at $J$ to $M_0+S_0L$ at $K$. 

This then is the geometric picture of the state of equilibrium: the forces and moments are represented by a triangular loop in 4D, and for the local moments, there is a triangle which grows in area along the bar. The book-keeping necessary to quantify all information may have been somewhat involved, but the final geometric picture of an evolving triangle is actually rather simple. 

\begin{figure} [ht!] \centering
\includegraphics[width = 0.6\textwidth]{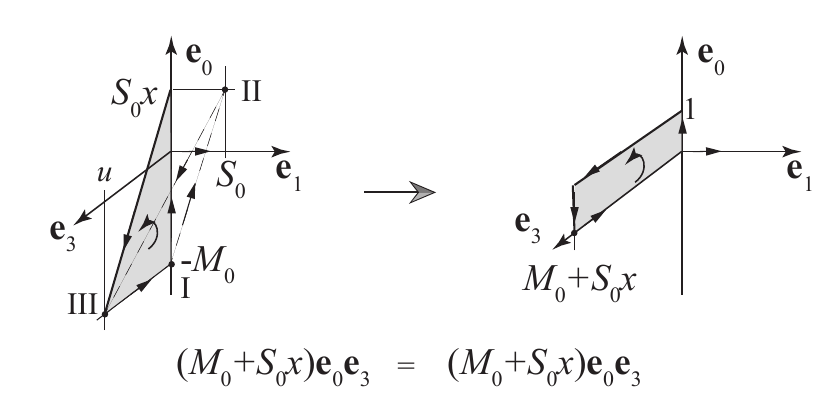}
\caption{The flag bivector $\be_0 \wedge (M_0 +S_0 x) \be_3$that will be used as shorthand for the bending moment.}
\label{MSx}
\end{figure}

As we proceed to introduce kinematics in the next section, we shall require the full 4D picture. In order to simplify matters and to keep the diagrams readily intelligible, we shall treat the bending moment at any point $x$ as the vector $(M_0 + S_0x) \be_3$, this being merely a convenient surrogate for the larger geometrical, triangular object projected onto the $\be_0\be_3$ plane. Essentially, the simple trapezoidal moment diagram of Fig.~\ref{basicbeamloop}d is used as shorthand for the full geometrical picture of Fig.~\ref{ForceDiagramforBeam}. The wedge product of $\be_0$ with that simpler diagram gives the same vector and bivector quantities as the fuller picture (see Fig.~\ref{MSx}). 

\subsubsection{The compatibility system}
Consider now the following compatibility system. Let the beam have a small lateral displacement $v(x)\be_2$, such that the deformed beam is a curve lying in the $xy$  plane with bases $\be_1$ and $\be_2$.

\begin{figure} [ht!] \centering
\includegraphics[width = 135mm]{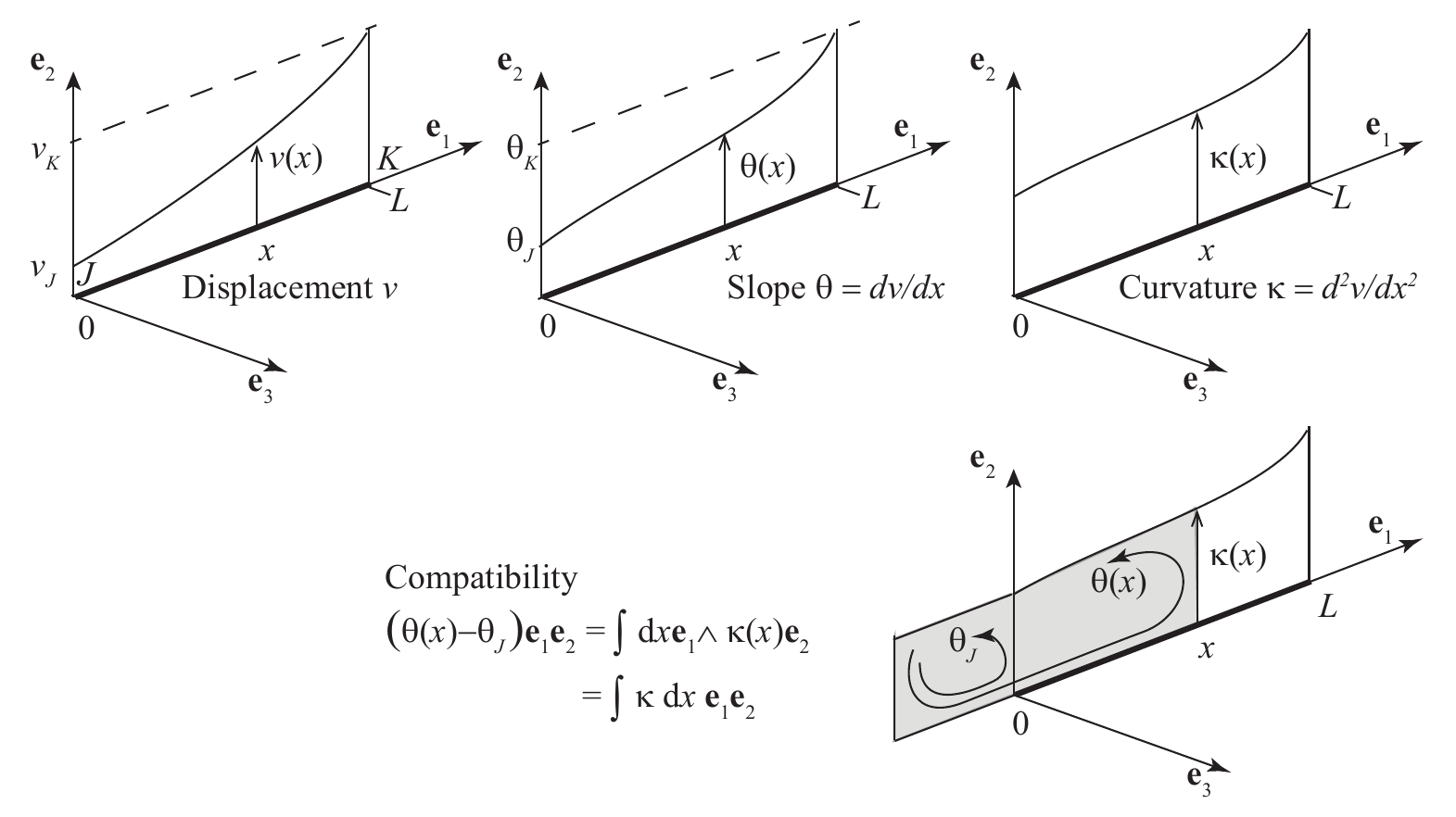}
\caption{The kinematic variables of displacement, slope and curvature shown as graphs along the beam. Slopes $\theta $ may also be represented as oriented areas under the curvature graph. The initial slope $\theta_I$ at end $I$  may be represented by an oriented area $\theta_I \be_1 \be_2$ to the left of the origin.}
\label{IntByParts3A}
\end{figure}

\begin{figure} [ht!] \centering
\includegraphics[width = 100mm]{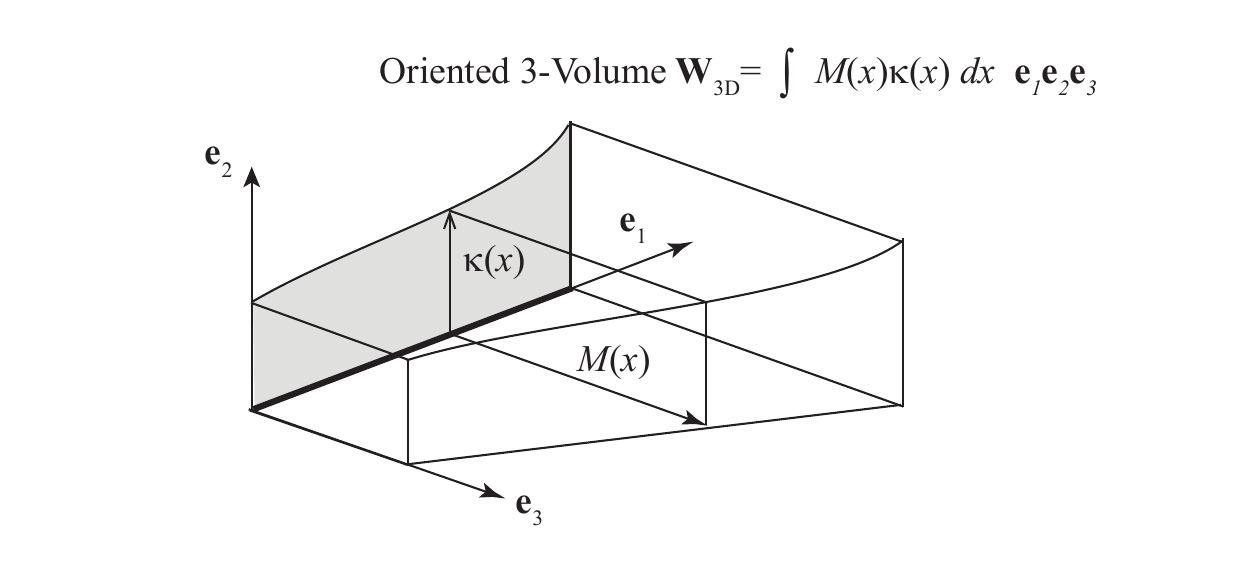}
\caption{The work integral $\int_0^L M(x) \kappa(x) dx $ represented as a 3-volume along the beam.}
\label{IntByParts3C0}
\end{figure}

We may plot the curvature $\kappa(x) = d^2v/dx^2$ as a graph on the $xy$ plane, as shown in Fig.~\ref{IntByParts3A}. Let the slope of the beam be  $dv/dx = \theta(x)$, and  let the beam have slope $\theta(0) = \theta_J$ at the left-hand node $J$ where $x=0$. We represent this initial slope geometrically by a rectangle of oriented area $\bTheta_J = \theta_J \be_1\be_2$ on the $xy$ plane, located just to the left of the origin. In Fig.~\ref{IntByParts3A}, this rectangle has been drawn so as to neatly adjoin the $\kappa(x)$ graph at $x=0$, even though such neatness is not necessary.

The slope $\theta(x)$ at any general point $x$ along the beam is then given by the integral of the graph of the curvature $\kappa = d^2v/dx^2$, plus the initial slope $\theta_J$.
The geometric representation as a bivector $\bTheta(x) = \theta(x) \be_1 \be_2$ is shown in Fig.~\ref{IntByParts3A}, this being the oriented area below the curvature graph, together with the adjoining initial rectangle.
The slope of the beam will thus be represented by an  $\be_1\be_2$ bivector.

The sign convention is defined in the statement
\begin{equation}
\bTheta(x) - \bTheta_J = \int_0^x \lbrace dx' \be_1 \wedge \kappa(x') \be_2 \rbrace = \int_0^x \kappa(x')dx' \ \be_1\be_2
\end{equation}

For the  equilibrium system, there is a bending moment $\beem(x) = M(x) \be_3 = (M_0 +S_0 x) \be_3$ about the $z$ axis. As explained in the previous section, in the full loop representation this will be represented by a bivector on the $\be_0 \be_3$ plane, but we use the vector in direction $\be_3$ for shorthand here (Fig.~\ref{MSx}).

\subsubsection{The Virtual Work}
Fig.~\ref{IntByParts3C0} shows how the moment-curvature integral may now be constructed by summing thin rectangular slices having side lengths $M(x)$ by $\kappa(x)$ and thickness $dx$.  The integral $\int M \kappa dx$ along the beam is thus given by the oriented volume in 3-space shown in Fig.~\ref{IntByParts3C0}. We may define the orientation of this 3-volume representation of the Virtual Work as

\begin{equation}
\bW_{3D} =   \int_0^L \lbrace dx \be_1 \wedge \kappa(x) \be_2  \wedge M(x) \be_3 \rbrace =  \int_0^L M(x) \kappa(x) dx \ \be_1\be_2\be_3
\end{equation}

This 3-volume forms part of the geometrical realisation of the Virtual Work due to bending in the beam. Strictly, though, a dimension is missing. The moment $M(x)$ has been shown here as the vector $M\be_3$ as a shorthand for the more fundamental object, the bivector component $M\be_0\be_3$ of the loop in 4D.  If we represent this component by the simple flag bivector $\be_0 \wedge M\be_3$, then the 3-volume integral may be pre-multiplied by $\be_0$ to give the full representation of the Virtual Work as an oriented 4-volume:
\begin{equation}
\bW_{4D}  =   \int_0^L \lbrace dx \be_1 \wedge \kappa(x) \be_2  \wedge \left[ \be_0 \wedge M(x) \be_3 \right] \rbrace = \int_0^L M(x) \kappa(x) dx \ \be_0 \be_1 \be_2\be_3
 \end{equation}

\begin{figure} [ht!] \centering
\includegraphics[width = 120mm]{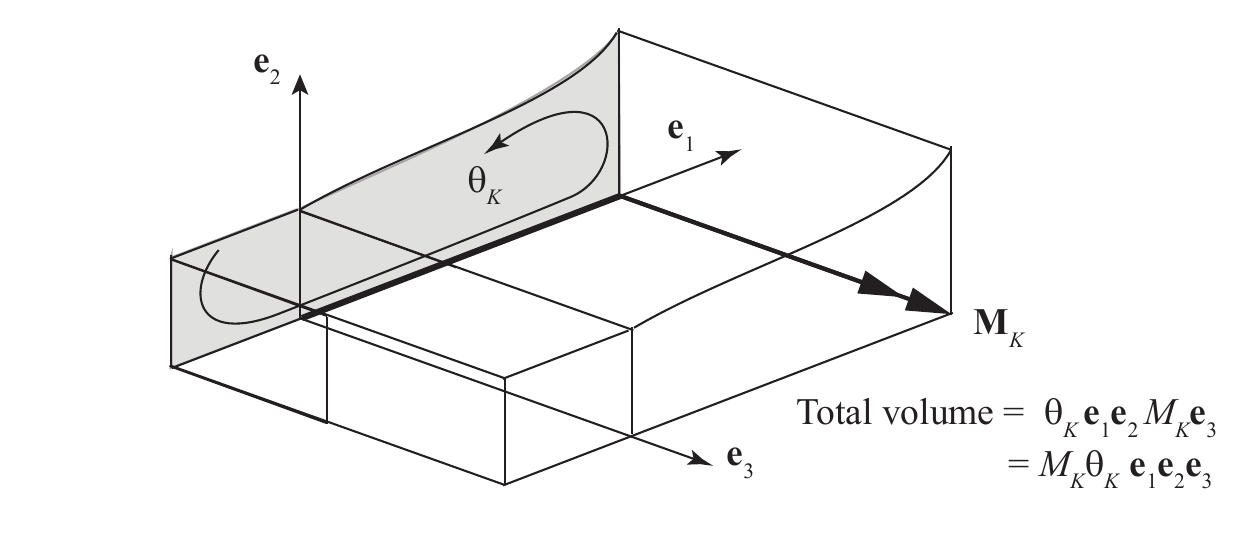}\\
\includegraphics[width = 120mm]{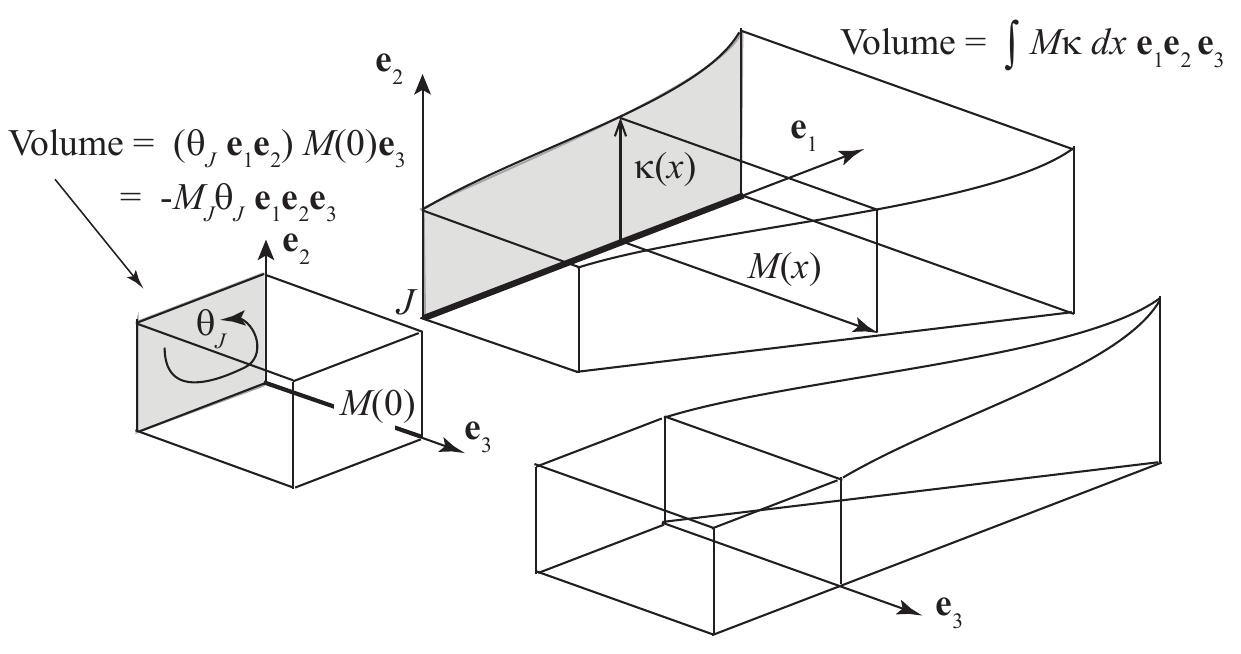}
\caption{Geometric decomposition of the work integrals. The overall volume represents the external virtual work
$M_K \theta_K $ due to moment and rotation at the right-hand end of the beam. The lower figures show this decomposed into three subvolumes:
one represents the external virtual work $M_J \theta_J $ due to moment and rotation at the left hand end of the beam.
Another is the internal virtual work due to flexure, given by the moment curvature integral  $\int_0^L M(x) \kappa(x) dx $ along the beam. The third volume represents the external virtual work due to the shear forces
and lateral displacements at the beam ends, as will be demonstrated. }
\label{IntByParts3C3D}
\end{figure}

For the present, we retain the 3-volume description of the Virtual Work. Fig.~\ref{IntByParts3C3D} shows how the outermost 3-volume $M_K\theta_K \be_1 \be_2 \be_3$ corresponding to the external virtual work involving the end moment and rotation at node $K$ may be decomposed into three sub-volumes. One of these has oriented volume $-M_J\theta_J \be_1 \be_2 \be_3$ associated with the initial (external) end moment and rotation at node $J$, one is associated with the (internal) moment curvature integral
$\int M \kappa dx$ and the third can be shown to involve the interaction of shear and lateral displacement.

\begin{figure} [ht!] \centering
\includegraphics[width = 120mm]{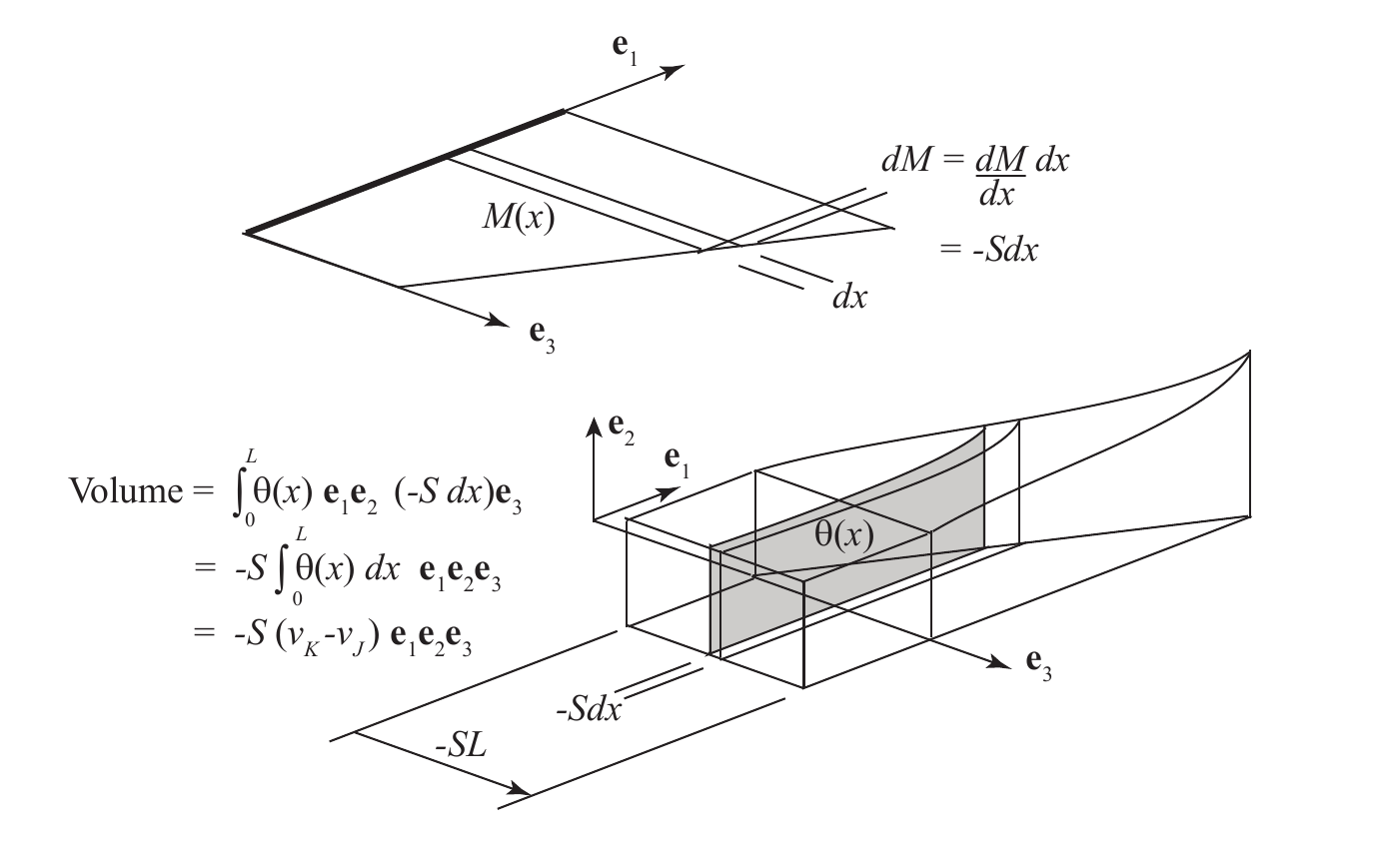}
\caption{The third sub-volume represents the external virtual work due to the shear forces
and lateral displacements at the beam ends. This follows by taking thin slices parallel to the $\be_1 \be_2$ plane.
These slices have area $\theta(x)$ and thickness $dM = -Sdx$. Since no forces are applied within the span of the beam, the shear force is constant. The volume integral is thus $-S\int \theta dx = -S (v_K - v_J)$ where $v_J$ and $v_K$ are the lateral displacements in the $\be_2$ direction at the ends $J$ and $K$.}
\label{IntByParts3E}
\end{figure}

Fig.~\ref{IntByParts3E} illustrates that the third sub-volume may be determined by summing slices of thickness $-Sdx$ and of area $\theta(x)$, such that the oriented volume is  $-S(v_K - v_J) \be_1\be_2 \be_3$.

\begin{figure} [b!] \centering
\includegraphics[width = 130mm]{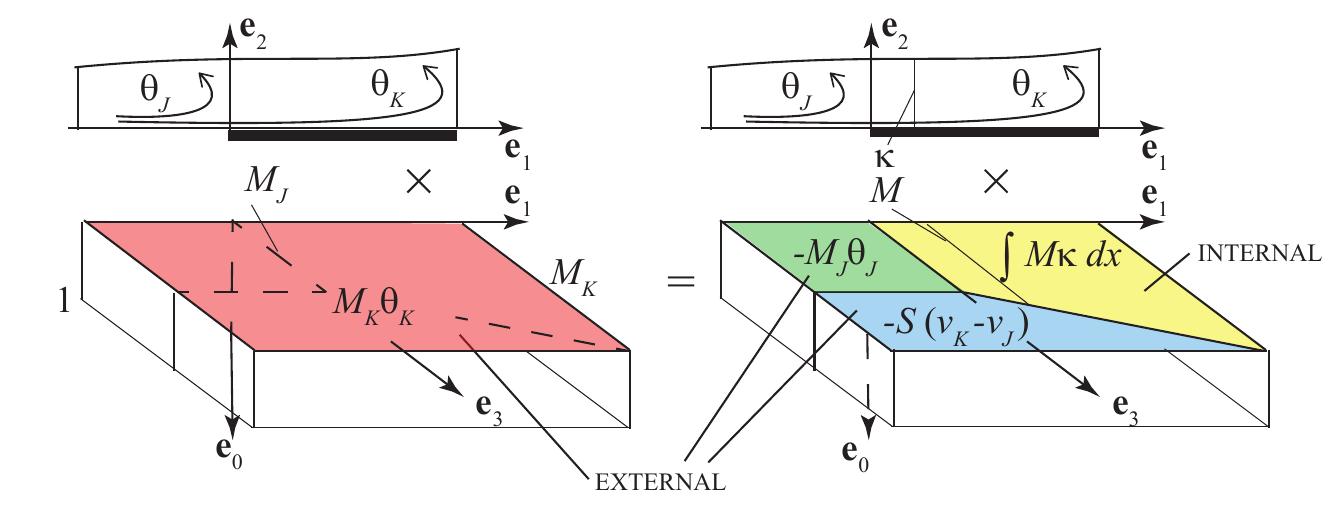}
\caption{A geometric statement of Virtual Work for beam flexure, expressed as the equivalence of 4D volumes. The beam is initially directed along the $\be_1$ axis, and flexes into the $\be_2$ direction.  Key to the compatibility system is the curvature graph, drawn in the $\be_2$ direction, whose integrals give rotation bivectors of the form $\theta \be_1\be_2$. Key to the equilibrium system are the moments, expressed as vectors of magnitude $M$ in the $\be_3$ direction, but pre-multiplied by the unit vector $\be_0$ to create moment bivectors of the form $M\be_0\be_3$. The various element of the Virtual Work are colour-coded, with Red, Green and Blue for external contributions and Yellow for the Internal Virtual Work. Care must be taken with signs. Here, Internal (Yellow) equals External (Red minus Green minus Blue).}
\label{IntByParts3F2}
\end{figure}

Equating the overall volume (in its 4-volume form) with the sum of its three sub-volumes, we thus obtain
\begin{equation}
M_K\theta_K \be_0 \be_1 \be_2 \be_3   = \left( -M_J\theta_J \ + \ \int_0^L M(x) \kappa(x) dx \ - \ S (v_K - v_J)  \right) \be_0\be_1 \be_2\be_3
\end{equation}
which rearranges to
\begin{equation}
\left( M_K\theta_K  + M_J\theta_J  + S_K v_K  + S_J v_J \right) \be_0 \be_1 \be_2 \be_3   =  \int_0^L M(x) \kappa(x) dx \   \be_0\be_1 \be_2\be_3
\end{equation}
where the left-hand side is associated with the external actions and the right hand side is  associated with the internal actions. This is of course the familiar statement of Virtual Work for a beam loaded at its ends, but here we give geometric expression to the statement that is usually only written algebraically.

It is, of course, also merely a geometric statement of integration by parts
\begin{eqnarray}
\int_0^L M(x) \kappa(x) dx  \ = \ \int_0^L \ M \frac{d \theta}{dx} dx  &  = & \left[ M \theta \right]_0^L - \int_0^L  \frac{dM}{dx} \theta dx \nonumber \\ 
& = &  \left[ M \theta \right]_0^L +  S \int_0^L  \theta dx \nonumber \\
& = &  M_K\theta_K  + M_J\theta_J  + S_K v_K  + S_J v_J
\end{eqnarray}
This states that internal virtual work must equal external virtual work, and there is no requirement that the equilibrium system $M(x)$ causes the curvatures $\kappa(x)$ of the compatibility system.

The geometry of the Principle of Virtual Work for 2D beam flexure is encapsulated in Fig.~\ref{IntByParts3F2}. The volume associated with the internal virtual work integral is coloured yellow and the various positive and negative contributions to the external virtual work are coloured red, green and blue. The $\be_0\be_1\be_2\be_3$ 4-volumes are created by the wedge products of the $\be_0\be_3$ bivectors representing moment with the $\be_1\be_2$ bivectors representing rotation.

\subsection{Canonical Examples of Flexure (Jan 2019)}

\begin{figure} [h!] \centering
\includegraphics[width = 130mm]{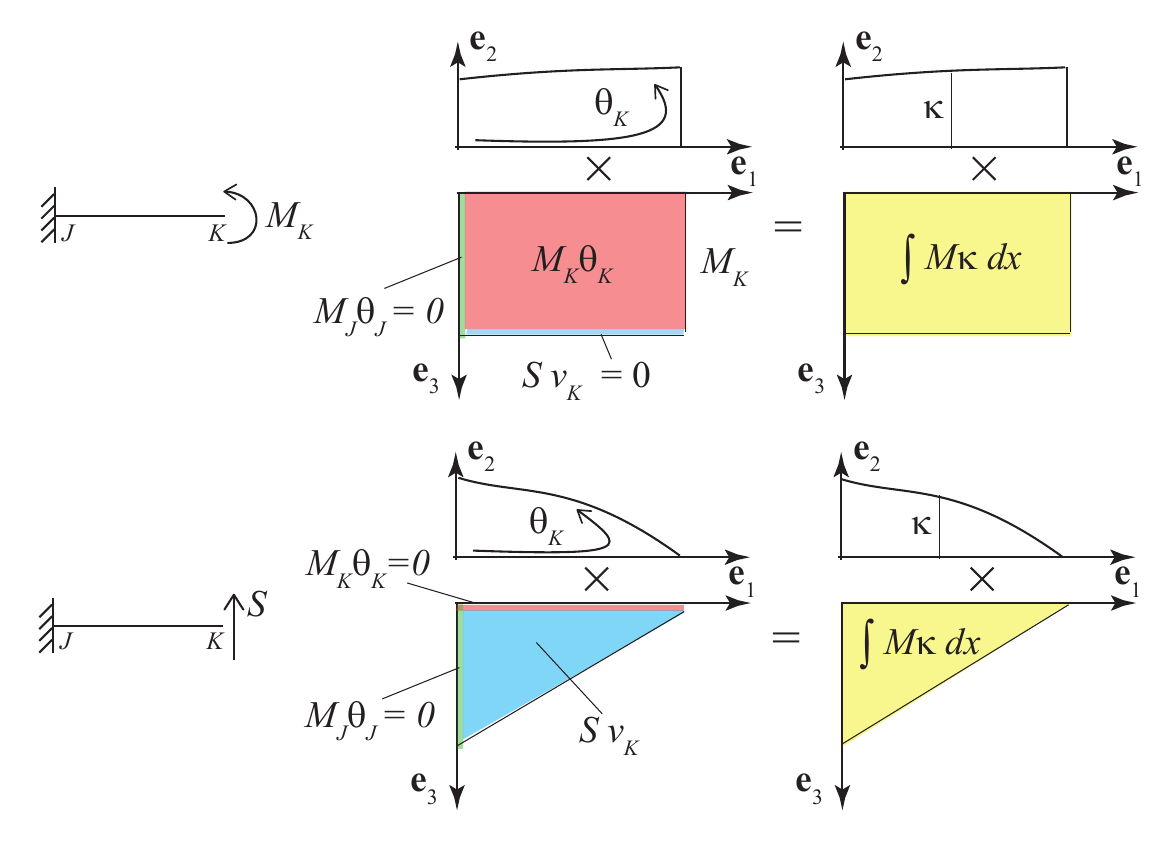}
\caption{The geometry of Virtual Work behind two canonical configurations of beam flexure - a cantilever with end moment and with tip load.}
\label{TwoExamples}
\end{figure}

Fig.~\ref{TwoExamples} shows the geometry of Virtual Work associated with two canonical configurations of beam flexure - a cantilever with end moment and with tip load. In each case, the partition of the Virtual Work into the various external and internal volumes is illustrated. At this stage, as in the previous section, no particular Material Law has been adopted, and the curvature graphs that define the compatibility systems are completely general.

\begin{figure} [ht!] \centering
\includegraphics[width = 130mm]{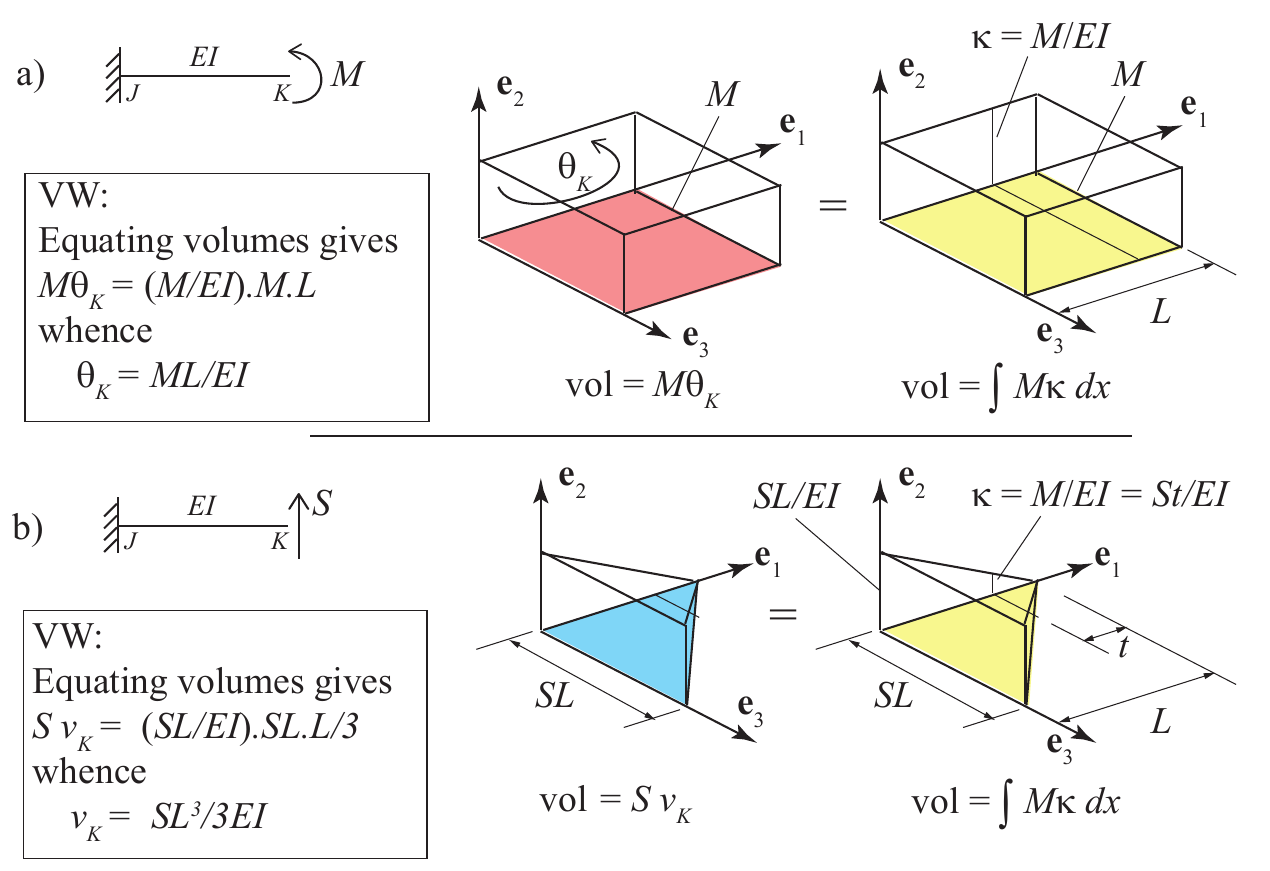}
\caption{Two familiar canonical configurations of flexure of a uniform elastic beam. a) shows a cantilever with end moment and b) shows a cantilever carrying a transverse end load. In each case, the Force Method equates the 4-volumes of internal and external Virtual Work, with the equilibrium and compatibility systems being the same.   
In each case, the (trivial) $\be_0$ contribution to the moment bivector has not been shown. }
\label{TwoExamples3D}
\end{figure}

In Fig.~\ref{TwoExamples3D}, an elastic material law has been assumed, with flexural rigidity $EI$ constant along the beam. The compatibility and equilibrium systems are now linked, with the curvature given by $\kappa = M/EI$. In each case, the volume integrals are evident. For the applied end moment (Fig.~\ref{TwoExamples3D}a), there is a cuboid of Internal Virtual Work having side lengths $(M/EI) \times  M \times L$. This must equal the volume $M_K\theta_K$ of the External Virtual Work cuboid, leading immediately to the familiar result that $\theta_K = M_KL/EI$.

For the cantilever with tip load (Fig.~\ref{TwoExamples3D}b), the volumes are pyramidal. The Internal Virtual Work is thus one third of the enclosing cuboid of sides $(SL/EI) \times  SL \times L$. This
must equal the External Virtual Work volume $Sv_K$, leading to the tip deflection $v_K = SL^3/EI$.

In these diagrams, the factor $\be_0$ in the moment bivector has been omitted for brevity. We thus obtain a reduced (3D) geometrical description which is easier to visualise, but lacks the more general power of the full 4D description, where the moment is represented by some more general oriented bivector.

\subsection{The $\Gamma$ frame}
\begin{figure} [ht!] \centering
\includegraphics[width = 150mm]{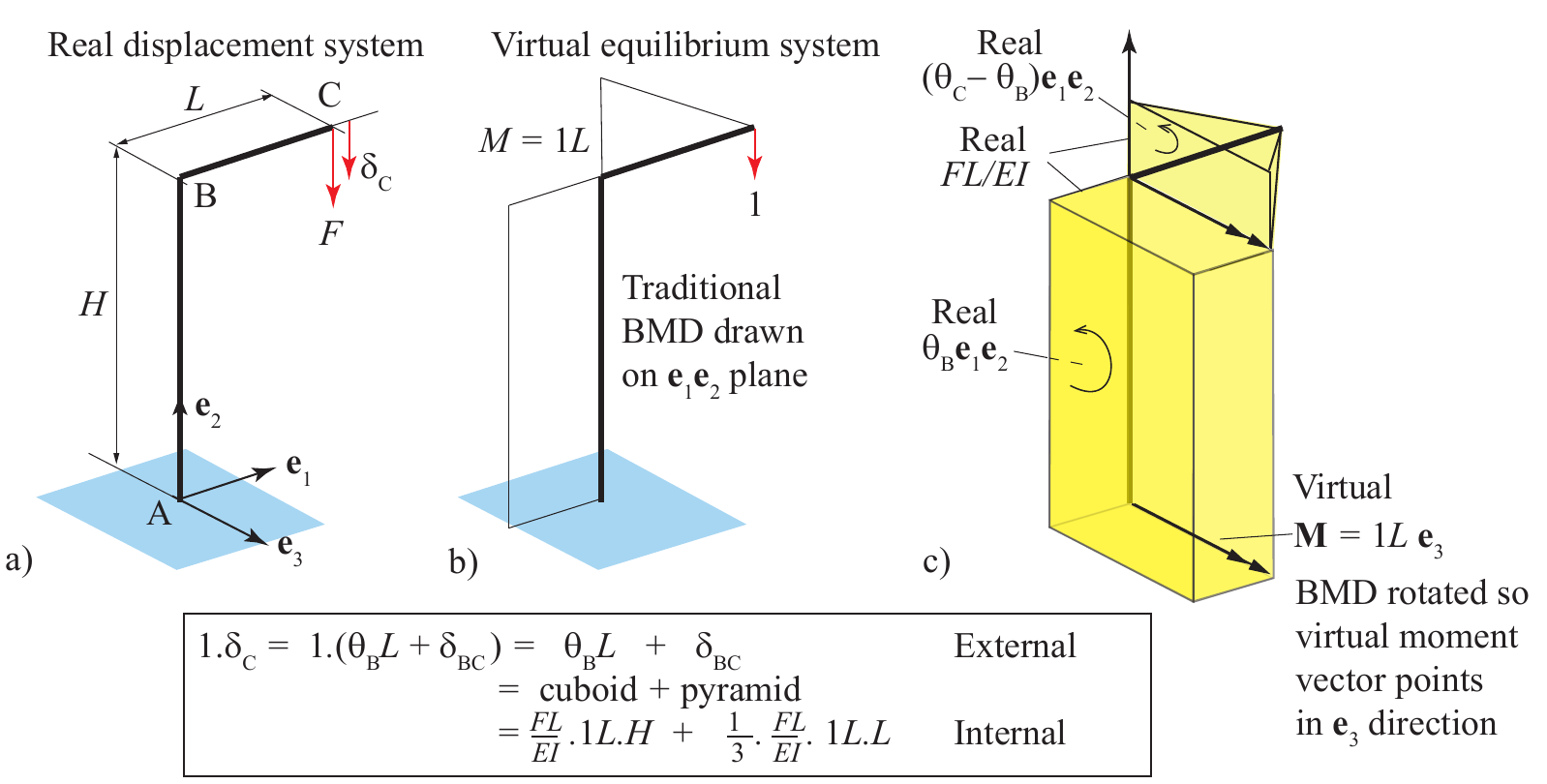}
\caption{The geometry underlying the Virtual Work calculation giving the tip deflection of a loaded $\Gamma$ frame. 
a), b) show the real displacement and virtual equilibrium systems respectively. c) shows the volume integrals which lead to the result. The unit vector in the direction $\be_0$ of the loop on $\be_0\be_3$ representing the virtual moment is not shown. 
}
\label{Gamma}
\end{figure}

A loaded $\Gamma$ frame, perhaps a rudimentary crane, provides a simple example of flexure (see Fig.~\ref{Gamma}a). A typical question would be, given the vertical load $F$ at the tip, find the vertical component $\delta$ of the deflection there. Following the Force Method, as usual, a virtual system is envisaged  which has a unit vertical load to pluck out the displacement component at the tip. 
The Virtual Work integrals are shown in Fig.~\ref{Gamma}c.  Graphs of the real curvatures are plotted on the $\be_1\be_2$ plane. The traditional Bending Moment Diagram for virtual moments is shown in  Fig.~\ref{Gamma}b, but for the volume integrals of Fig.~\ref{Gamma}c, this has been rotated so that the virtual moments plot as vectors in the $\be_3$ direction. This is a notational shorthand for the virtual moment loop on the $\be_0\be_3$ plane. As in Fig.~\ref{TwoExamples3D} earlier, the external and internal work volumes coincide, and equating them leads to  
\begin{equation*}
\delta_C =  \frac{FL^2}{EI}\left(H + \frac{L}{3}\right)
\end{equation*}

This is the usual result as would be obtained by more traditional methods, but here we have added a graphical representation of the various contributions to the Virtual Work $\bW = W\be_0\be_1\be_2\be_3$ as 4-volumes in the full space.

\section{Summary and Conclusions}
It has been demonstrated how the loop formalism  previously applied to describe the equilibrium of forces and moments in 3D rigid-jointed frame structures~\cite{McRobieArXiv1,McRobieArXiv2} could be applied in an analogous manner to describe their displacements and rotations. The structural state is then given by three sets of loops in 4D: form loops, force loops and displacement loops. 
Although it may at first seem unusual and perhaps even awkward to describe a rotation as the projected bivector area of a loop in 4D, it was shown how this may correspond to the common practice of integrating graphs of beam curvatures to obtain changes in end rotations. For beam structures with a linear elastic material law, the graph of the Bending Moment Diagram may be interpreted as a loop, and beam curvatures are scaled versions of this.  The Principle of Virtual Work then manifests itself as an equivalence of 4-volumes given by the wedge product of the force loop with the displacement loop.

\bibliography{McRobieBib3}

\end{document}